\documentclass[12pt]{article}
\usepackage{amsmath, epsfig, cite}
\usepackage{amssymb}
\usepackage{amsfonts}
\usepackage{latexsym}
\usepackage{graphicx}
\usepackage{url}
\newtheorem{thm}{Theorem}[section]
\newtheorem{prop}[thm]{Proposition}

\newtheorem{defi}[thm]{Definition}
\newtheorem{lem}[thm]{Lemma}

\newcommand{\qed}{{\hfill\rule{4pt}{7pt}}}
\def\pf{\noindent {\it Proof.} }
\numberwithin{equation}{section}

\makeatletter \@addtoreset{equation}{section} \makeatother
\begin{document}
\rule{0cm}{1cm}
\begin{center}
{\Large\bf Ramanujan-type Congruences for

Broken $2$-Diamond Partitions Modulo $3$}
\end{center}

 \vskip 2mm \centerline{William Y.C. Chen$^1$,  Anna R.B. Fan$^2$
 and Rebecca T. Yu$^3$ }

\begin{center}
$^{1,2,3}$Center for Combinatorics, LPMC-TJKLC\\
Nankai University, Tianjin 300071, P. R. China

$^{1}$Center for Applied Mathematics\\
 Tianjin University, Tianjin 300072, P. R. China

\vskip 2mm
 Email: $^1$chen@nankai.edu.cn, $^2$fanruice@mail.nankai.edu.cn,
 $^3$yuting\_shuxue@mail.nankai.edu.cn
\end{center}

{\noindent \bf Abstract.} { The notion of broken $k$-diamond partitions was
introduced by Andrews and Paule. Let  $\Delta_k(n)$ denote the number of
broken $k$-diamond partitions of $n$. They also posed three conjectures on
the congruences of $\Delta_2(n)$ modulo $2$, $5$ and $25$. Hirschhorn and
sellers proved the conjectures for modulo $2$,
 and Chan proved cases of modulo $5$.
 For the case of modulo 3, Radu and Sellers obtained an
  infinite family of congruences for $\Delta_2(n)$.
  In this paper, we obtain two infinite families of
 congruences for $\Delta_2(n)$ modulo $3$ based on a
 formula of Radu and Sellers, the 3-dissection formula of
 the generating function of triangular number due to Berndt, and
 the properties of the $U$-operator, the $V$-operator, the Hecke operator and
 the Hecke eigenform. For example, we find that $\Delta_2(243n+142)\equiv \Delta_2(243n+223)\equiv0\pmod{3}$. The infinite family of Radu and Sellers
 and the  two infinite families derived
 in this paper have two congruences in common, namely, $\Delta_2(27n+16)\equiv\Delta_2(27n+25)\equiv0 \pmod{3}$.

 \vskip 3mm

 \vskip 3mm


\section{Introduction}

The objective of this paper is to derive two infinite families of
congruences for the number of broken $2$-diamond partitions modulo 3. The
notion of broken $k$-diamond partitions was   introduced by Andrews and
Paule \cite{G.AndrewsP.Paule}. A broken $k$-diamond partition
$\pi=(a_1,a_2,a_3,\ldots;b_2,b_3,b_4,\ldots)$ is a plane partition
satisfying the  relations illustrated in Figure \ref{fig1}, where $a_i, b_i$
are nonnegative integers and $a_i\rightarrow a_j$ means $a_i\geq a_j$. More
precisely,  each building block in Figure \ref{fig1}, except for the broken
block $(b_2,b_3,\ldots b_{2k+2})$  has the same order structure as shown in
Figure \ref{fig3}. We call each block a $k$-elongated partition diamond of
length $1$, or a $k$-elongated diamond, for short.

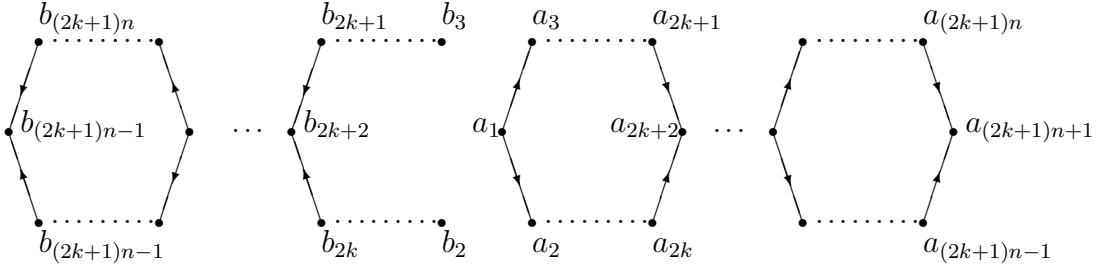
\begin{figure}[h]
\begin{center}
\setlength{\unitlength}{0.8mm}
\begin{picture}(120,60)
\put(60,30){$a_1$}
\put(65,30){\circle*{1.5}}
\put(65,30){\vector(1,3){3}}
\put(68,39){\line(1,3){2}}
\put(70,48){$a_3$}
\put(70,45){\circle*{1.5}}
\multiput(72,45)(2,0){9}{.}
\put(90,48){$a_{2k+1}$}
\put(90,45){\circle*{1.5}}
\put(90,45){\vector(1,-3){2.8}}
\put(92.8,36.6){\line(1,-3){2.2}}
\put(83,30){$a_{2k+2}$}
\put(95,30){\circle*{1.5}}
\put(95,30){\line(-1,-3){2.2}}
\put(90,15){\vector(1,3){2.8}}
\put(90,10){$a_{2k}$}
\put(90,15){\circle*{1.5}}
\multiput(88,15)(-2,0){9}{.}
\put(70,10){$a_2$}
\put(70,15){\circle*{1.5}}
\put(70,15){\line(-1,3){2}}
\put(65,30){\vector(1,-3){3}}
\multiput(100,30)(2,0){3}{.}
\put(110,30){\circle*{1.5}}
\put(110,30){\vector(1,3){3}}
\put(113,39){\line(1,3){2}}
\put(115,45){\circle*{1.5}}
\multiput(117,45)(2,0){9}{.}
\put(135,48){$a_{(2k+1)n}$}
\put(135,45){\circle*{1.5}}
\put(135,45){\vector(1,-3){2.8}}
\put(137.8,36.6){\line(1,-3){2.2}}
\put(142,30){$a_{(2k+1)n+1}$}
\put(140,30){\circle*{1.5}}
\put(140,30){\line(-1,-3){2.2}}
\put(135,15){\vector(1,3){2.8}}
\put(135,10){$a_{(2k+1)n-1}$}
\put(135,15){\circle*{1.5}}
\multiput(133,15)(-2,0){9}{.}
\put(115,15){\circle*{1.5}}
\put(115,15){\line(-1,3){2}}
\put(110,30){\vector(1,-3){3}}
\put(32,30){$b_{2k+2}$}
\put(30,30){\circle*{1.5}}
\put(30,30){\line(1,3){2}}
\put(35,45){\vector(-1,-3){3}}
\put(35,48){$b_{2k+1}$}
\put(35,45){\circle*{1.5}}
\multiput(37,45)(2,0){9}{.}
\put(55,48){$b_3$}
\put(55,45){\circle*{1.5}}
\put(55,10){$b_2$}
\put(55,15){\circle*{1.5}}
\multiput(53,15)(-2,0){9}{.}
\put(35,10){$b_{2k}$}
\put(35,15){\circle*{1.5}}
\put(35,15){\vector(-1,3){3}}
\put(30,30){\line(1,-3){2}}
\multiput(24,30)(-2,0){3}{.}
\put(-15,30){$b_{(2k+1)n-1}$}
\put(-17,30){\circle*{1.5}}
\put(-17,30){\line(1,3){2}}
\put(-12,45){\vector(-1,-3){3}}
\put(-12,48){$b_{(2k+1)n}$}
\put(-12,45){\circle*{1.5}}
\multiput(-10,45)(2,0){9}{.}
\put(8,45){\circle*{1.5}}
\put(8,45){\line(1,-3){2.2}}
\put(13,30){\vector(-1,3){2.8}}
\put(13,30){\circle*{1.5}}
\put(13,30){\vector(-1,-3){2.8}}
\put(8,15){\line(1,3){2.2}}
\put(8,15){\circle*{1.5}}
\multiput(6,15)(-2,0){9}{.}
\put(-12,10){$b_{(2k+1)n-1}$}
\put(-12,15){\circle*{1.5}}
\put(-12,15){\vector(-1,3){3}}
\put(-17,30){\line(1,-3){2}}
\end{picture}
\caption{A broken $k$-diamond of length $2n$.}\label{fig1}
\end{center}
\end{figure}

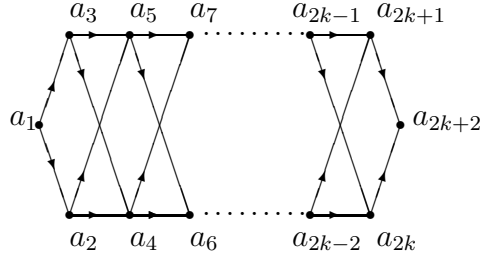
\begin{figure}[h]
\begin{center}
\setlength{\unitlength}{0.8mm}
\begin{picture}(100,50)
\put(10,30){$a_1$}
\put(15,30){\circle*{1.5}}
\put(15,30){\vector(1,3){2.5}}
\put(17.5,37.5){\line(1,3){2.5}}
\put(20,48){$a_3$}
\put(20,45){\circle*{1.5}}
\put(20,45){\vector(1,0){5}}
\put(25,45){\line(1,0){5}}
\put(20,45){\vector(1,-3){2.5}}
\put(22.5,37.5){\line(1,-3){7.5}}
\put(30,48){$a_5$}
\put(30,45){\circle*{1.5}}
\put(30,45){\vector(1,0){5}}
\put(35,45){\line(1,0){5}}
\put(30,45){\vector(1,-3){2.5}}
\put(32.5,37.5){\line(1,-3){7.5}}
\put(40,48){$a_7$}
\put(40,45){\circle*{1.5}}
\multiput(42,45)(2,0){9}{.}
\put(57,48){$a_{2k-1}$}
\put(60,45){\circle*{1.5}}
\put(60,45){\vector(1,0){5}}
\put(65,45){\line(1,0){5}}
\put(60,45){\vector(1,-3){2.5}}
\put(62.5,37.5){\line(1,-3){7.5}}
\put(71,48){$a_{2k+1}$}
\put(70,45){\circle*{1.5}}
\put(70,45){\vector(1,-3){2.5}}
\put(72.5,37.5){\line(1,-3){2.5}}
\put(77,30){$a_{2k+2}$}
\put(75,30){\circle*{1.5}}
\put(75,30){\line(-1,-3){2.5}}
\put(70,15){\vector(1,3){2.5}}
\put(71,10){$a_{2k}$}
\put(70,15){\circle*{1.5}}
\put(60,15){\vector(1,0){5}}
\put(65,15){\line(1,0){5}}
\put(57,10){$a_{2k-2}$}
\put(60,15){\circle*{1.5}}
\put(60,15){\vector(1,3){2.5}}
\put(62.5,22.5){\line(1,3){7.5}}
\multiput(58,15)(-2,0){9}{.}
\put(40,10){$a_6$}
\put(40,15){\circle*{1.5}}
\put(30,15){\vector(1,0){5}}
\put(35,15){\line(1,0){5}}
\put(30,10){$a_4$}
\put(30,15){\circle*{1.5}}
\put(30,15){\vector(1,0){5}}
\put(25,15){\line(1,0){5}}
\put(20,15){\vector(1,0){5}}
\put(30,15){\vector(1,3){2.5}}
\put(32.5,22.5){\line(1,3){7.5}}
\put(20,10){$a_2$}
\put(20,15){\circle*{1.5}}
\put(20,15){\line(-1,3){2.5}}
\put(15,30){\vector(1,-3){2.5}}
\put(20,15){\vector(1,3){2.5}}
\put(22.5,22.5){\line(1,3){7.5}}
\end{picture}
\caption{A $k$-elongated   diamond.}\label{fig3}
\end{center}
\end{figure}

For example, Figure \ref{fig4} gives a broken $2$-diamond partition \[ \pi=(10,8,9,7,6,3,2,1,0,1;3,5,2,1,1,1).\]

\begin{figure}[h]
\begin{center}
\setlength{\unitlength}{0.8mm}
\begin{picture}(100,50)
\put(48,29){$10$}
\put(55,30){\circle*{1}}
\put(55,30){\vector(1,2){3}}
\put(58,36){\line(1,2){3}}
\put(61,43){$9$}
\put(61,42){\circle*{1}}
\put(61,42){\vector(1,0){6}}
\put(67,42){\line(1,0){6}}
\put(61,42){\vector(1,-2){4}}
\put(65,34){\line(1,-2){8}}
\put(73,43){$6$}
\put(73,42){\circle*{1}}
\put(73,42){\vector(1,-2){3}}
\put(76,36){\line(1,-2){3}}
\put(55,30){\vector(1,-2){3}}
\put(58,24){\line(1,-2){3}}
\put(61,13){$8$}
\put(61,18){\circle*{1}}
\put(61,18){\vector(1,0){6}}
\put(67,18){\line(1,0){6}}
\put(61,18){\vector(1,2){4}}
\put(65,26){\line(1,2){8}}
\put(73,13){$7$}
\put(73,18){\circle*{1}}
\put(73,18){\vector(1,2){3}}
\put(76,24){\line(1,2){3}}
\put(81,29){$3$}
\put(79,30){\circle*{1}}

\put(79,30){\vector(1,2){3}}
\put(82,36){\line(1,2){3}}
\put(85,43){$1$}
\put(85,42){\circle*{1}}
\put(85,42){\vector(1,0){6}}
\put(91,42){\line(1,0){6}}
\put(85,42){\vector(1,-2){4}}
\put(89,34){\line(1,-2){8}}
\put(97,43){$1$}
\put(97,42){\circle*{1}}
\put(79,30){\vector(1,-2){3}}
\put(82,24){\line(1,-2){3}}
\put(85,13){$2$}
\put(85,18){\circle*{1}}
\put(85,18){\vector(1,0){6}}
\put(91,18){\line(1,0){6}}
\put(85,18){\vector(1,2){4}}
\put(89,26){\line(1,2){8}}
\put(97,13){$0$}
\put(97,18){\circle*{1}}
\put(25,29){$1$}
\put(23,30){\circle*{1}}
\put(23,30){\circle*{1}}
\put(23,30){\line(1,2){3}}
\put(29,42){\vector(-1,-2){3}}
\put(29,43){$1$}
\put(29,42){\circle*{1}}
\put(29,42){\line(1,0){6}}
\put(41,42){\vector(-1,0){6}}
\put(29,42){\line(1,-2){8}}
\put(41,18){\vector(-1,2){4}}
\put(41,43){$5$}
\put(41,42){\circle*{1}}
\put(23,30){\line(1,-2){3}}
\put(29,18){\vector(-1,2){3}}
\put(29,13){$2$}
\put(29,18){\circle*{1}}
\put(29,18){\line(1,0){6}}
\put(41,18){\vector(-1,0){6}}
\put(29,18){\line(1,2){8}}
\put(41,42){\vector(-1,-2){4}}
\put(41,13){$3$}
\put(41,18){\circle*{1}}
\put(23,30){\vector(-1,-2){3}}
\put(20,24){\line(-1,-2){3}}
\put(17,13){$1$}
\put(17,18){\circle*{1}}

\end{picture}
\caption{A broken 2-diamond partition of 60.}\label{fig4}
\end{center}
\end{figure}
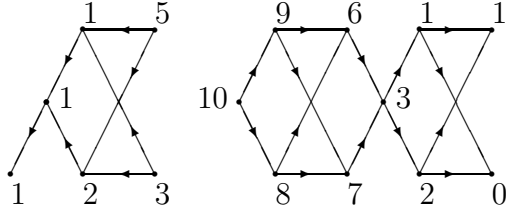

  Let $\Delta_k(n)$ denote the number of broken $k$-diamond partitions
  of $n$, and let $B_k(q)$ denote the generating function for $\Delta_k(n)$, that is
  \[ B_k(q)=\sum_{n\geq0}\Delta_k(n)q^n.\]
 Andrews and Paule \cite{G.AndrewsP.Paule}
  obtained the following formula:
  \begin{align}\label{generate-delta} B_k(q) =\frac{(-q;q)_{\infty}}{(q;q)^2_{\infty}(-q^{2k+1};q^{2k+1})_{\infty}}.
\end{align}
Note that the above formula can also be written in terms of eta-quotients related
to modular forms:
\[B_k(q)=q^{(k+1)/12}\frac{\eta(2z)\eta((2k+1)z)}{\eta(z)^3\eta((4k+2)z)},\]
where $q=e^{2\pi iz}$.

From \eqref{generate-delta},  Andrews and Paule proved that for $n\geq 3$,
\begin{equation}
\Delta_1(2n+1)\equiv 0\pmod{3}.\label{delta1-3}
\end{equation}
They also posed three conjectures on congruences for the number of
 broken $2$-diamond partitions of $n$.
Fu \cite{shishuofu} and Mortenson \cite{EMortenson} found
 combinatorial proofs of  congruence \eqref{delta1-3}.
Meanwhile, Hirschhorn and Sellers \cite{MHirschhornSellers} gave a proof of
 \eqref{delta1-3} by deriving the following generating formula for
   $\Delta_1(2n+1)$
 \[\sum_{n\geq 0}\Delta_1(2n+1)q^n=3\frac{(q^2;q^2)_{\infty}^2 (q^6;q^6)_{\infty}^2}{(q;q)_{\infty}^6}.\]
Hirschhorn and Sellers also obtained the following congruences modulo 2
\begin{align}
&\Delta_1(4n+2)\equiv \Delta_1(4n+3)\equiv0 \pmod{2},\\[5pt]
&\Delta_2(10n+2)\equiv \Delta_2(10n+6)\equiv0 \pmod{2},\label{conj1}
\end{align}
where $n\geq 1$. The congruences in \eqref{conj1} were conjectured
by Andrews and Paule \cite{G.AndrewsP.Paule}.

 Chan \cite{shchan} found the following two infinite families of congruences modulo $5$ for broken
 $2$-diamond partitions:
\begin{align}
\Delta_2\left(5^{l+1}+\frac{3}{4}(5^l-1)+2\cdot5^l+1\right)\equiv  0\pmod{5},\label{chan1}\\[6pt]
\Delta_2\left(5^{l+1}+\frac{3}{4}(5^l-1)+4\cdot5^l+1\right) \equiv  0\pmod{5},\label{chan2}
\end{align}
where $l\geq1$.
 The cases  for $l=1$ in \eqref{chan1} and \eqref{chan2} were
conjectured by Andrews and Paule \cite{G.AndrewsP.Paule}, namely,
\begin{equation*}
\Delta_2(25n+14)\equiv\Delta_2(25n+24)\equiv0 \pmod{5}.
\end{equation*}

 Paule and Radu \cite{pauleradu} obtained the following infinite family
 of congruences modulo 5 for broken
$2$-diamond partitions. They showed that for any prime $p$ with $p\equiv
13,17 \pmod{20}$ and any nonnegative integer $n$,
 \begin{equation}\label{paule-radu-congruence}
\Delta_2\left(5 p^2n+4p-\frac{1}{4}(p-1)\right)\equiv0\pmod{5}.
\end{equation}
Moreover, they posed four conjectures on congruences for broken $3$-diamond
partitions and broken $5$-diamond partitions, which have been
confirmed by Xiong \cite{XXiong} and Jameson \cite{MJamson}.

For the broken $2$-diamond partitions, Radu and Sellers\cite{radusellers} showed
that
\begin{equation}\label{radu-1}
\sum_{n\geq0}\Delta_2(3n+1)q^n\equiv 2q\prod_{n\geq1}\frac{(1-q^{10n})^4}{(1-q^{5n})^2}\pmod{3},
\end{equation}
which implies the following congruence relations
\begin{align}
&\Delta_2(15n+1)\equiv\Delta_2(15n+7)\equiv0 \pmod{3},\label{15-1}\\[5pt]
&\Delta_2(15n+10)\equiv\Delta_2(15n+13)\equiv0 \pmod{3},\label{15-2}
\end{align}
and
\begin{equation}\label{radu-seller-congruence}
\Delta_2\left(3p^2n+\frac{3}{4}\left(p(4k+3)-1\right)+1\right)\equiv 0\pmod{3},
\end{equation}
where prime $p\equiv3\pmod{4}$, $0\leq k \leq p-1$ and $k\neq
\frac{p-3}{4}$.

In this paper, we use \eqref{radu-1} to establish two new infinite families of congruences of $\Delta_2(n)$ modulo $3$ by using the 3-dissection
 formula of the generating function of triangular numbers and properties of  the $U$-operator, the $V$-operator, the Hecke operator
  and the Hecke eigenform.
\begin{thm}\label{thm}
For $l\geq1$, we have
\begin{align}
&\Delta_2\left(3^{2l+1}n+\frac{3}{4}(3^{2l}-1)+3^{2l}+1\right)\equiv 0 \pmod{3},\label{thmeq1}\\[6pt]
&\Delta_2\left(3^{2l+1} n+\frac{3}{4}(3^{2l}-1)+2\cdot3^{2l}+1\right)\equiv 0\pmod{3}.\label{thmeq2}
\end{align}
\end{thm}

\section{Preliminaries}

In this section, we give an overview of some definitions and properties of
modular forms which will be used in the proof of Theorem \ref{thm}.
For more details on the theory of modular forms, see for example, Ono \cite{Omodular}.
In this paper, we shall be concerned with the modular forms for the congruence subgroup $\Gamma_0(N)$ of $SL_2(\mathbb{Z})$,
 where
 \[\Gamma_0(N):=\left\{\left.\left(
                                         \begin{array}{cc}
                                           a & b \\
                                           c & d \\
                                         \end{array}\right)\in SL_2(\mathbb{Z})
                                       \right|c\equiv0\pmod{N}
\right\}.\]
First, recall the definition of a modular form.

\begin{defi}
Let $k$ be a positive integer and $\chi$ be a Dirichlet character modulo
$N$, we say that a meromorphic function $f(z):\mathbb{H}\rightarrow
\mathbb{C}$ is a modular form for $\Gamma_0(N)$ with Nebentypus $\chi$ and
weight $k$, if
\begin{description}
\item[\rm (1)]$f(z)$ is a holomorphic function on $\mathbb{H}$.
\item[\rm (2)]$f(z)$ satisfies the following relation for $z \in\mathbb{H}$ and $\gamma\in \Gamma_0(N)$
\[f(\gamma z)=\chi(d)(cz+d)^kf(z).\]
\item[\rm (3)]$f(z)$ is holomorphic at the cusps of $\Gamma_0(N)$,
here the cusps are the equivalent classes of $Q\cup\{\infty\}$
on the action of $\Gamma_0(N)$.
\end{description}
\end{defi}

The set of modular forms for
$\Gamma_0(N)$   with Nebentypus $\chi$ and weight $k$ is denoted by $M_k(\Gamma_0(N),\chi)$.
When $\chi$ is trivial, that is, $\chi(d)=1$ for all integers $d$, we denote $M_k(\Gamma_0(N),\chi)$ by $M_k(\Gamma_0(N))$.
In addition, if $f(z)$ vanishes at the cusps of $\Gamma_0(N)$, then we say that
$f(z)$ is a cusp form.
We denote by $S_k(\Gamma_0(N),\chi)$ the set of cusp forms for $\Gamma_0(N)$ with Nebentypus $\chi$ and weight $k$.

To prove Theorem \ref{thm}, we need the operators $U$, $V$ and the Hecke
operator.

\begin{defi}\label{defi-operation}
Let $ f(z)=\sum_{n\geq0}a(n)q^n$ be a modular form in
$M_k(\Gamma_0(N),\chi)$. If $d$ is a positive integer, then the
\emph{$U$-operator} $U(d)$ and the \emph{$V$-operator} $V(d)$ are defined by
\begin{equation}\label{defi-U-operator}
f(z)| U(d) :=\sum_{n\ge 0} a(d n) q^n
\end{equation}
and
\begin{equation}\label{defi-V-operator}
f(z)|V(d) := \sum_{n\ge 0} a(n) q^{d n}.
\end{equation}
For any positive integer $m$, the action of the $m$-th Hecke operator
$T_{m,k,\chi}:=T(m)$ is given by
\begin{align*}
f(z)\mid T(m):=\sum_{n\geq 0}\left(\sum_{d|\gcd(m,n)}\chi(d)d^{k-1}a(mn/d^2) \right)q^n.
\end{align*}
In particular, for any prime $p$
\begin{align}\label{mf12}
f(z)\mid T(p):=\sum_{n\geq 0}\left( a(pn)+\chi(p)p^{k-1}a\left(n/p\right) \right)q^n.
\end{align}
\end{defi}
From the definition of the $U$-operator and the Hecke operator, one sees that
for any prime $p$, we have
\begin{equation}
f(z)|U(p)\equiv f(z)|T(p) \pmod{p}.\label{up-tp}
\end{equation}

The operators $U$, $V$ and the Hecke operator
have the following properties.

\begin{prop}\label{prop-operation}
Suppose that $f(z)\in M_k(\Gamma_0(N),\chi)$.
\begin{description}
\item[\rm(1)]If $d$ is a positive integer, then
\[f(z)|V(d)\in M_k(\Gamma_0(Nd),\chi).\]
\item[\rm(2)]If $d|N$, then
\[f(z)|U(d)\in M_k(\Gamma_0(N),\chi).\]
\item[\rm(3)]If $m$ is a positive integer, then
\[f(z)|T(m)\in M_k(\Gamma_0(N),\chi).\]
\end{description}
\end{prop}

A Hecke eigenform associated with the Hecke operators is defined as follows.

\begin{defi}\label{defi-eigenform}
A modular form $f(z)\in M_k(\Gamma_0(N),\chi)$ is called a Hecke eigenform
if for every $m\geq 2$ there is a complex number $\lambda(m)$ for which
\[f(z)|T(m)=\lambda(m)f(z).\]
\end{defi}

If a Hecke eigenform   is a cusp form, then the following proposition can be
used to compute the value $\lambda(m)$.

\begin{prop}\label{prop-eigenform}
Suppose that $f(z)=\sum_{n\geq0}a(n)q^n$ is a cusp form in
$S_k(\Gamma_0(N),\chi)$ with $a(1)=1$. If $f(z)$ is an eigenform, then for $m\geq 1$,
\[f(z)|T(m)=a(m)f(z).\]

\end{prop}

\section{Proof of the Main Theorem}

In this section, we give a proof of Theorem \ref{thm} by using the
 approach of Chan \cite{shchan}. We first establish the following congruence on the
generating function of $\Delta_2(3n+1)$.
Let $\psi(q)$ be the generating function for the triangular numbers, that is,
\begin{equation}
\psi(q)=\sum_{n\geq 0}q^{n(n+1)/2}=\prod_{n\geq1}\frac{(1-q^{2n})^2}{(1-q^{n})}.\label{psi}
\end{equation}

\begin{lem}\label{lem1}
We have
\begin{align}\label{u3}
\psi(q^{15})^2\sum_{n\geq0}\Delta_2(3n+1)q^{n+4}\equiv 2q^5\psi(q^5)^8\pmod 3.
 \end{align}
\end{lem}

\noindent {\it{Proof.}}
By \eqref{radu-1} and  \eqref{psi}, we see that
\begin{align}\label{propeq}
\sum_{n\geq0}\Delta_2(3n+1)q^n\equiv 2q\prod_{n\geq1}\frac{(1-q^{10n})^4}{(1-q^{5n})^2}=2q\psi(q^5)^2\pmod{3}.
\end{align}
It follows that
\begin{align}
 2q^5\psi(q^5)^8&=q^4\psi(q^5)^6\cdot(2q\psi(q^5)^2)\notag\\[6pt]
  &\equiv\psi(q^5)^6\sum_{n\geq0}\Delta_2(3n+1)q^{n+4} \pmod 3.\label{psi1}
\end{align}
Since
\begin{equation*}
\frac{(q;q)_{\infty}^3}{(q^3;q^3)_\infty}\equiv 1 \pmod{3},
\end{equation*}
see, Ono \cite{Omodular},  we find that
\begin{align}
\psi(q^5)^6\sum_{n\geq0}\Delta_2(3n+1)q^{n+4}
 \equiv\psi(q^{15})^2\sum_{n\geq0}\Delta_2(3n+1)q^{n+4}\pmod 3.\label{psi2}
 \end{align}
Combining \eqref{psi1} and \eqref{psi2}, we obtain that
\begin{align*}
\psi(q^{15})^2\sum_{n\geq0}\Delta_2(3n+1)q^{n+4}\equiv 2q^5\psi(q^5)^8\pmod 3,
 \end{align*}
as desired.\qed

It is not difficult to show that $q^5\psi(q^5)^8$ is a Hecke eigenform.

\begin{lem}\label{eigen-lemma}
The function $q^5\psi(q^5)^8$ is a Hecke eigenform in
$M_4(\Gamma_0(10))$. More precisely, for $m\geq 2$, we have
\[
q^5\psi(q^5)^8|T(m)=\lambda(m)q^5\psi(q^5)^8,
\]
where $\lambda(m)=[q^m]q\psi(q)^8$, that is, the coefficient of $q^m$ in
   $q\psi(q)^8$.
\end{lem}

\pf
It has been shown by  Chan \cite{H.H.Chan} that $q\psi(q)^8$ is a Hecke eigenform in
$M_4(\Gamma_0(2))$.
 By Definition \ref{defi-eigenform},   for $m\geq 2$, there exits a complex number $\lambda(m)$ such that
 \begin{equation}
 q\psi(q)^8|T(m)=\lambda(m)q\psi(q)^8. \label{hecke-equation}
 \end{equation}
 Observing that $q\psi(q)^8$ is a cusp form for which the coefficient
 of $q$ is $1$,
 by Proposition \ref{prop-eigenform}  we find
 \begin{equation}\label{lambda(m)}
 \lambda(m)=[q^m]q\psi(q)^8.
 \end{equation}
Substituting $q$ by $q^5$ in \eqref{hecke-equation}, we obtain that
\[q^5\psi(q^5)^8|T(m)=\lambda(m)q^5\psi(q^5)^8.\]
 Meanwhile, by Proposition \ref{prop-operation} (2), we have
 $$q^5\psi(q^5)^8=q\psi(q)^8|V(5)\in M_4(\Gamma_0(10)).$$
Thus $q^5\psi(q^5)^8$ is a Hecke eigenform in
$M_4(\Gamma_0(10))$. \qed

 We are now ready to prove Theorem \ref{thm}.

\noindent {\it{Proof of Theorem \ref{thm}.}}
Let $f(z)=2q^5\psi(q^5)^8$ with $q=e^{2\pi i z}$.
Applying the $U$-operator $U(3)$  to $f(z)$ gives
 \begin{align*}
 2q^5\psi(q^5)^8|U(3)&\equiv \left.\left(\psi(q^{15})^2\sum_{n\geq0}
 \Delta_2(3n+1)q^{n+4}\right)\right|U(3)\pmod 3\notag\\[6pt]
 &=\psi(q^5)^2\sum_{n\geq0}\Delta_2(9n-11)q^{n}\hskip 21mm\pmod 3.
 \end{align*}
 On the other hand, Lemma \ref{eigen-lemma} implies that
 \begin{align*}
 q^5\psi(q^5)^8|T(3)&=\left([q^3]q\psi(q)^8\right)q^5\psi(q^5)^8 \nonumber\\[6pt]
 &=28q^5\psi(q^5)^8\equiv q^5\psi(q^5)^8 \pmod 3.
 \end{align*}
Employing relation  \eqref{up-tp}, we deduce that
\begin{align*}
\psi(q^5)^2\sum_{n\geq0}\Delta_2(9n-11)q^{n}&\equiv 2q^5\psi(q^5)^8|U(3)&\pmod 3\\
&\equiv 2q^5\psi(q^5)^8|T(3)\equiv 2q^5\psi(q^5)^8&\pmod 3.
\end{align*}
Substituting $n$ by $n+2$ in the above congruence, we get
\begin{equation}\label{9generate}
\sum_{n\geq0}\Delta_2(9n+7)q^{n}\equiv 2q^3\psi(q^5)^6\equiv 2q^3\psi(q^{15})^2\pmod 3.
\end{equation}
By the $3$-dissection formula
 of $\psi(q)$ due to Berndt \cite[p49]{B.C.Berndt}, one sees that
 $\psi(q)$ can be written in the following form
\begin{equation}\label{3-dissection}
\psi(q)=A(q^3)+q\psi(q^9),
\end{equation}
where $A(q)$ is power series  in $q$.
Plugging \eqref{3-dissection} into the right hand side of \eqref{9generate},
we find
\begin{align}
\sum_{n\geq0}\Delta_2(9n+7)q^{n}&\equiv 2q^3\left(A(q^{45})+q^{15}\psi(q^{135})\right)^2 \hskip 38mm \pmod{3}\notag\\
&\equiv 2q^3\left(A(q^{45})^2+2q^{15}A(q^{45})\psi(q^{135})+q^{30}\psi(q^{135})^2\right) \pmod{3}.\label{9generate-trans}
\end{align}
Extracting those terms whose powers of $q$ are congruent to $6$ modulo $9$ in \eqref{9generate-trans}, we have
\begin{equation}\label{9n+6}
\sum_{n\geq0}\Delta_2\left(9(9n+6)+7\right)q^{9n+6}\equiv 2q^{33}\psi(q^{135})^2 \pmod{3}.
\end{equation}
Dividing both sides of \eqref{9n+6} by $q^6$ and replacing $q^9$ by $q$,
we deduce that
\begin{equation}\label{81generate}
\sum_{n\geq0}\Delta_2(3^4n+61)q^{n}\equiv 2q^{3}\psi(q^{15})^2 \pmod{3}.
\end{equation}
Combining \eqref{9generate} and \eqref{81generate}, we arrive at
\begin{equation}\label{repeat}
\sum_{n\geq0}\Delta_2(3^2n+7)q^{n}\equiv\sum_{n\geq0}\Delta_2(3^4n+61)q^{n}\pmod{3}.
\end{equation}
Iterating \eqref{repeat} with $n$ replaced by $9n+6$, we conclude that
for $l\geq 1$,
\begin{equation*}
 \sum_{n\geq0}\Delta_2(3^2n+7)q^{n}\equiv\sum_{n\geq 0}\Delta_2\left(3^{2l}n+\frac{3}{4}(3^{2l}-1)+1\right)q^{n}\pmod{3}.
\end{equation*}
From (\ref{9generate}) together with the above relation, it follows that
\begin{align}\label{final}
 \sum_{n\geq 0}\Delta_2\left(3^{2l}n+\frac{3}{4}(3^{2l}-1)+1\right)q^{n}\equiv 2q^{3}\psi(q^{15})^2\pmod{3}.
\end{align}
Since there are no terms with powers of $q$   congruent
to 1, 2 modulo 3 in $2q^{3}\psi(q^{15})^2$,
substituting $n$ by $3n+1$ and $3n+2$ respectively in
\eqref{final}, we obtain the following  infinite families of congruences
\begin{align*}
\Delta_2\left(3^{2l+1}n+\frac{3}{4}(3^{2l}-1)+3^{2l}+1\right)&\equiv0\pmod{3},\\[6pt]
\Delta_2\left(3^{2l+1}n+\frac{3}{4}(3^{2l}-1)+2\cdot 3^{2l}+1\right)&\equiv 0 \pmod{3}.
\end{align*}
This completes the proof.\qed

Here are some examples of Theorem \ref{thm}.
For $n\geq 0$, we have
\begin{align}
&\Delta_2(27n+16)\equiv \Delta_2(27n+25)\equiv0\pmod{3},\label{27}\\[5pt]
&\Delta_2(243n+142)\equiv \Delta_2(243n+223)\equiv0\pmod{3},\\[5pt]
&\Delta_2(2187n+1276)\equiv \Delta_2(2187n+2005)\equiv0\pmod{3}.
\end{align}
Notice that  congruences in \eqref{27} are also contained in the
infinite family of congruences \eqref{radu-seller-congruence} derived by Radu and Sellers.

\vspace{0.5cm}
 \noindent{\bf Acknowledgments.}  This work was supported by  the 973
Project, the PCSIRT Project of the Ministry of Education,  and the
National Science Foundation of China.

\end{document}